\documentclass[11pt,reqno]{amsart}
\usepackage{amssymb,verbatim,enumerate,ifthen}
\usepackage[mathscr]{eucal}
\usepackage[utf8]{inputenc}
\usepackage[T1]{fontenc}
\allowdisplaybreaks[1]
\def\N{\mathbb{N}}
\def\R{\mathbb{R}}

\newtheorem{theorem}{Theorem}
\newtheorem*{theorem*}{Theorem}
\def\Thm#1#2{\ifthenelse{\equal{#1}{*}}{\begin{theorem*}#2\end{theorem*}}
             {\begin{theorem}\label{T#1}#2\end{theorem}}}
\newtheorem{Atheorem}{Theorem}

\def\thm#1{Theorem~\ref{T#1}}

\newtheorem{proposition}[theorem]{Proposition}
\newtheorem*{proposition*}{Proposition}
\def\Prp#1#2{\ifthenelse{\equal{#1}{*}}{\begin{proposition*}#2\end{proposition*}}
             {\begin{proposition}\label{P#1}#2\end{proposition}}}

\newtheorem{corollary}[theorem]{Corollary}
\newtheorem*{corollary*}{Corollary}
\def\Cor#1#2{\ifthenelse{\equal{#1}{*}}{\begin{corollary*}#2\end{corollary*}}
             {\begin{corollary}\label{C#1}#2\end{corollary}}}
\def\cor#1{Corollary~\ref{C#1}}

\newtheorem{lemma}[theorem]{Lemma}
\newtheorem*{lemma*}{Lemma}
\def\Lem#1#2{\ifthenelse{\equal{#1}{*}}{\begin{lemma*}#2\end{lemma*}}
             {\begin{lemma}\label{L#1}#2\end{lemma}}}
\def\lem#1{Lemma~\ref{L#1}}
\newtheorem{Alemma}{Lemma}

\theoremstyle{definition}
\newtheorem{remark}[theorem]{Remark}
\newtheorem*{remark*}{Remark}
\def\Rem#1#2{\ifthenelse{\equal{#1}{*}}{\begin{remark*}\rm #2\end{remark*}}
             {\begin{remark}\label{R#1}\rm #2\end{remark}}}

\newtheorem{example}[theorem]{Example}
\newtheorem*{example*}{Example}
\def\Exa#1#2{\ifthenelse{\equal{#1}{*}}{\begin{example*}\rm #2\end{example*}}
             {\begin{example}\label{Ex#1}\rm #2\end{example}}}

\def\eq#1{{\rm(\ref{E#1})}}
\def\Eq#1#2{\ifthenelse{\equal{#1}{*}}
  {\begin{equation*}\begin{aligned}#2\end{aligned}\end{equation*}}
  {\begin{equation}\begin{aligned}\label{E#1}#2\end{aligned}\end{equation}}}

\def\sign{\mathop{\hbox{\rm sign}}\nolimits}

\def\transp#1{\mathop{\hbox{$[\,#1\,]^T$}}\nolimits}

\def\diag{\mathop{\hbox{\rm diag}}\nolimits}

\def\comment#1{}

\begin{document}
\vspace{5mm}

\date{\today}

\title[Inequalities for nonsymmetric generalized Bajraktarevi\'c means]{Local and global Hölder- and Minkowski-type inequalities \\ for nonsymmetric generalized Bajraktarevi\'c means}

\author[R.\ Gr\"unwald]{Rich\'ard Gr\"unwald}
\address[R.\ Gr\"unwald]{Institute of Mathematics and Computer Science, University of Nyíregyháza, H-4400 Nyí\-regyháza, Pf. 166, Hungary and Doctoral School of Mathematical and Computational Sciences, University of Debrecen, H-4002 Debrecen, Pf. 400}
\email{grunwald.richard@nye.hu, richard.grunwald@science.unideb.hu}

\author[Zs. P\'ales]{Zsolt P\'ales}
\address[Zs. P\'ales]{Institute of Mathematics, University of Debrecen, H-4002 Debrecen, Pf. 400, Hungary}
\email{pales@science.unideb.hu}

\thanks{The research of the first author was supported by the \'UNKP-22-3 New National Excellence Program of the Ministry for Innovation and Technology from the source of the National Research, Development and Innovation Fund and by the Scientific Council of the University of Nyíregyháza. The research of the second author was supported by the K-134191 NKFIH Grant.}


\begin{abstract}
The aim of this paper is to investigate inequalities that are analogous to the Minkowski and Hölder inequalities by replacing the addition and the multiplication by a more general operation, and instead of using power means, generalized Bajraktarević means are considered, in particular, Gini means. A further aim is to introduce the concept of local and global validity of such inequalities and to characterize them in both senses.
\end{abstract}

\maketitle

\section{Introduction}

Throughout this paper, the symbols $\N,\R,$ and $\R_+$ will stand for the sets of natural (i.e., positive integer), real, and positive real numbers, respectively, and $I$ will always denote a nonempty open real interval. Let $n,k\in\N$. In the sequel, the $i$th entry of a real vector
\Eq{*}{
	x:=(x_i)_{i\in\{1,\dots,n\}}=
	\begin{pmatrix}
		x_1 \\
		\vdots \\
		x_n
	\end{pmatrix}\in\R^n
}
will be denoted by $x_i$, and analogously, 
the $i$th row and $j$th column of a real matrix
\Eq{*}{
	x:=(x_i^j)_{(i,j)\in\{1,\dots,n\}\times\{1,\dots,k\}}=
	\begin{pmatrix}
		x_1^1 & \cdots & x_1^k \\
		\vdots & & \vdots \\
		x_n^1 & \cdots & x_n^k
	\end{pmatrix}\in\R^{n\times k}
}
will be denoted by $x_i$ and $x^j$, respectively. For convenience, we identify $\R^{n\times k}$ by $(\R^n)^k$ in the standard manner. We define the \emph{transpose} $[x]^T\in\R^{k\times n}$ of the matrix $x\in\R^{n\times k}$ by 
\Eq{*}{
	[x]^T:=(x_j^i)_{(i,j)\in\{1,\dots,n\}\times\{1,\dots,k\}}=
	\begin{pmatrix}
		x_1^1 & \cdots & x_n^1 \\
		\vdots & & \vdots \\
		x_1^k & \cdots & x_n^k.
	\end{pmatrix}
}
More generally, for a subset of matrices $X\subseteq\R^{n\times k}$,
the transpose $X^T$ of $X$ denotes the set $\{x^T\mid x\in X\}$.
Finally, the \emph{diagonal} $\diag(I^n)$ of $I^n$ is defined by 
\Eq{*}{
	\diag(I^n):=\{(x,\dots,x)\in\R^n\mid x\in I\}.
}
Let us introduce, for $n\in\N$, the \emph{diagonal map} $\Delta_n\colon\R\to\diag(\R^n)$ by
\Eq{*}{
	\Delta_n(x):=(x,\dots,x)\in\R^n
}
and, for an $n$-variable function $G\colon I^n\to\R$, the function $G^\Delta\colon I\to\R$ by
\Eq{*}{
	G^\Delta(x):=G(\Delta_n(x))\qquad(x\in I).
}
More generally, for $n,k\in\N$, we can define the map $\Delta_n^k\colon\R^k\to(\diag(\R^n))^k\subseteq\R^{n\times k}$: if $y\in\R^k$, then let $\Delta_n^k(y)$ denote the $n\times k$ matrix whose $j$th column equals $\Delta_n(y_j)$ for $j\in\{1,\dots,k\}$. Whenever a regularity property is assumed to be valid at each point of the domain of a function, then we do not emphasize the set on which the property in question holds. 

The celebrated inequalities discovered by Hölder and Minkowski can be formulated in various contexts, for instance, in the setting of power (or Hölder) means.

To recall the standard Hölder(--Rogers) inequality (which was discovered by Rogers in 1888 and by Hölder in 1889), let $p,q>1$ with $p^{-1}+q^{-1}=1$. Then, for all $n\in\N$ and $x,y\in\R_+^n$, the inequality
\Eq{*}{
  \frac{x_1y_1+\dots+x_ny_n}{n}
  \leq\bigg(\frac{x_1^p+\dots+x_n^p}{n}\bigg)^{\frac{1}{p}}
  \bigg(\frac{y_1^q+\dots+y_n^q}{n}\bigg)^{\frac{1}{q}}
}
is valid. In the particular case $p=q=2$, this inequality reduces to the so-called Cauchy--Bunyakovsky--Schwarz inequality, which in the above form was established by Cauchy in 1821. Given a real parameter $p\geq1$, the standard Minkowski inequality (established in 1910) states that the $p$th power mean is subadditive, i.e., for all $n\in\N$ and $x,y\in\R_+^n$, the inequality
\Eq{*}{
	\bigg(\frac{(x_1+y_1)^p+\dots+(x_n+y_n)^p}{n}\bigg)^{\frac{1}{p}}
	\leq\bigg(\frac{x_1^p+\dots+x_n^p}{n}\bigg)^{\frac{1}{p}}
	+\bigg(\frac{y_1^p+\dots+y_n^p}{n}\bigg)^{\frac{1}{p}}
}
holds.

Briefly, the aim of this paper is to investigate analogous inequalities by replacing the addition and the multiplication by more general operations, and instead of power means, also using generalized Bajraktarević means and, in particular, Gini means. A further aim is to introduce the concept of local and global validity of such inequalities and to characterize them in both senses.

Let $n\in\N$. Given a strictly monotone continuous function $f\colon I\to\R$ and an $n$-tuple of positive valued functions $p=(p_1,\dots,p_n)\colon I\to\R_{+}^n$, the \emph{$n$-variable nonsymmetric generalized Bajraktarević mean} $A_{f,p}\colon I^n\to I$ is given by the following formula:
\Eq{*}{
    A_{f,p}(x):=
	f^{-1}\bigg(\frac{p_1(x_1)f(x_1)+\dots+p_n(x_n)f(x_n)}{p_1(x_1)+\dots+p_n(x_n)}\bigg) \qquad (x\in I^n).
}
This is an extension of the notion introduced by Bajraktarevi\'c in the symmetric setting in \cite{Baj58} and \cite{Baj63}, that is, in the case when $p_1=\dots=p_n$, i.e., when all the \emph{weight functions} are the same. In the sequel, the sum of these weight functions will be denoted by $p_0$, i.e., $p_0:=p_1+\dots+p_n$. It is easy to see that $A_{f,p}$ is a strict mean, i.e.,
\Eq{*}{
	\min\{x_1,\dots,x_n\}\leq A_{f,p}(x)\leq \max\{x_1,\dots,x_n\}\qquad (x\in I^n)
}
holds, and the inequalities are strict if $\min\{x_1,\dots,x_n\}<\max\{x_1,\dots,x_n\}$. The equality and comparison problem of nonsymmetric generalized Bajraktarević means have been investigated by the authors in the recent papers \cite{GruPal20} and \cite{GruPal22a}. 

The main goal of this article is to investigate Hölder- and Minkowski-type inequality problems for the $n$-variable nonsymmetric generalized Bajraktarević means.
More generally, we are going to derive necessary as well as sufficient conditions for the local as well as for the global validity of the functional inequality
\Eq{Mineq}{
	M_0(\Phi(x_1),\dots,\Phi(x_n))
	\leq \Phi(M_1(x^1),\dots,M_k(x^k)),
}
where $n,k\in\N$, for $\alpha\in\{0,\dots,k\}$, $I_\alpha\subseteq\R$ is a nonempty open interval, $I:=I_1\times\cdots\times I_k$, $M_\alpha\colon I_\alpha^n\to I_\alpha$ is an $n$-variable mean and $\Phi\colon I\to I_0$. If there exists an open set $U\subseteq I^n$ such that $\diag(I^n)\subseteq U$ and \eq{Mineq} holds for all $x\in U^T\subseteq\prod_{\alpha=1}^k I_\alpha^n$, then we say that \eq{Mineq} holds in the local sense. If \eq{Mineq} is valid for all $x\in (I^n)^T=\prod_{\alpha=1}^k I_\alpha^n$, then we say that \eq{Mineq} holds in the global sense. Clearly, the global validity of \eq{Mineq} implies its local validity.

Then we consider the particular case of \eq{Mineq} when all the means are $n$-variable nonsymmetric generalized Bajraktarević means, i.e., we consider the inequality 
\Eq{ineq}{
	A_{f_0,p^0}(\Phi(x_1),\dots,\Phi(x_n))&\leq \Phi(A_{f_1,p^1}(x^1),\dots,A_{f_k,p^k}(x^k)),
}
where $n,k\in\N$, for $\alpha\in\{0,\dots,k\}$, $f_\alpha\colon I_\alpha\to\R$ is a strictly monotone continuous function, $p^\alpha\colon I_\alpha\to\R_+^n$. We obtain necessary as well as sufficient conditions for its validity in the local and also in the global sense.

We mention some important particular cases of \eq{ineq}.
\begin{enumerate}
 \item If $k=1$, $I_0=I_1=:I$ and $\Phi(x)=x$, then \eq{ineq} reduces to the local and global comparison problem of nonsymmetric generalized Bajraktarević means.
 \item If $k\geq2$, $I_0=I_1=\dots=I_k=:I$, $\Phi(x_1,\dots,x_k)=\frac1k(x_1+\dots+x_k)$, and $f_0=f_1=\dots=f_k=:f$,
 $p^0=p^1=\dots=p^k=:p$, then \eq{ineq} means the Jensen convexity of $A_{f,p}$. In this case, \eq{ineq} is said to be a Jensen-type inequality.
 \item If $k\geq2$, $I_0=I_1=\dots=I_k=\R_+$, $\Phi(x_1,\dots,x_k)=x_1+\dots+x_k$, and $f_0=f_1=\dots=f_k=:f$,
 $p^0=p^1=\dots=p^k=:p$, then \eq{ineq} expresses the subadditivity of $A_{f,p}$, which is often called a Minkowski-type inequality.
 \item If $k\geq2$, $I_0=I_1=\dots=I_k=\R_+$, $\Phi(x_1,\dots,x_k)=x_1\cdots x_k$, then \eq{ineq} reduces to a Hölder-type inequality for the means $A_{f_0,p^0},A_{f_1,p^1},\dots,A_{f_k,p^k}$.
\end{enumerate}

There are many results related to the Hölder- and Minkowski-type inequalities. Without completeness, we mention the following standard sources and the references therein: Hardy--Littlewood--Pólya \cite{HarLitPol34}, Beckenbach--Bellmann \cite{BecBel61},  Bullen--Mitrinović--Vasić \cite{BulMitVas88}, Mitrinović--Pečarić--Fink \cite{MitPecFin91}. We also quote the papers \cite{Los71a,Los71c,Los73a,Los77,Los83a,LosPal97a,LosPal11b} by Losonczi and the papers \cite{CziPal00,Pal82b,Pal83a,Pal83c,Pal86a}.

\section{Hölder- and Minkowski-type inequalities in the local sense}

For the investigation of inequality \eq{Mineq}, let us introduce the function $F\colon I_1^n\times\dots\times I_k^n \subseteq\R^{n\times k}\to\R$ by
\Eq{F}{
  F(x)=F(x^1,\dots,x^k)
  :=\Phi(M_1(x^1),\dots,M_k(x^k))-M_0(\Phi(x_1),\dots,\Phi(x_n)),
}
where $n,k\in\N$ and also set $I:=I_1\times\dots\times I_k$.

\Rem{F0}{
Observe that, for all $y\in I$, we have
\Eq{*}{
  F(\Delta_n^k(y))=0.
}
Indeed, by using the mean value property of $M_0,M_1,\dots,M_k$, it follows that
\Eq{*}{
  F(\Delta_n^k(y))&=F(\Delta_n(y_1),\dots,\Delta_n(y_k))
  =\Phi(M_1(\Delta_n(y_1)),\dots,M_k(\Delta_n(y_k)))-M_0(\Delta_n(\Phi(y)))\\
  &=\Phi(M_1^\Delta(y_1),\dots,M_k^\Delta(y_k))-M_0^\Delta(\Phi(y))=\Phi(y_1,\dots,y_k)-\Phi(y)=0.
}}

For the computation of the partial derivatives of $F$ at points of the form $\Delta_n^k(y)$, we formulate the following Lemma. In what follows, $\delta_{\cdot,\cdot}$ will stand for the standard Kronecker symbol.

\Lem{PDF}{Let $n,k\in\N$ with $n,k\geq 2$ and, for $\alpha\in\{0,\dots,k\}$, let $I_\alpha\subseteq\R$ be a nonempty open interval and $M_\alpha\colon I_\alpha^n\to I_\alpha$ be an $n$-variable mean, define $F\colon I_1^n\times\dots\times I_k^n\to\R$ by \eq{F} and let $\Phi\colon I\to I_0$.
\begin{enumerate}[(i)]
 \item Assume, for $\alpha\in\{0,\dots,k\}$, that $M_\alpha$ is partially differentiable on $\diag(I_\alpha^n)$ and that $\Phi$ is differentiable. Then, for all $i\in\{1,\dots,k\}$, $\ell\in\{1,\dots,n\}$, and $y\in I$,
\Eq{*}{
  \partial_{\ell+n(i-1)}F(\Delta_n^k(y))
  =\partial_i\Phi(y)\big(\partial_\ell M_i^\Delta(y_i)
  -\partial_\ell M_0^\Delta(\Phi(y))\big).
}
 \item Assume, for $\alpha\in\{0,\dots,k\}$, that $M_\alpha$ is twice partially differentiable on $\diag(I_\alpha^n)$ and that $\Phi$ is twice differentiable. Then, for all $i,j\in\{1,\dots,k\}$, $\ell,m\in\{1,\dots,n\}$, and $y\in I$,
\Eq{*}{
  \partial_{\ell+n(i-1)}\partial_{m+n(j-1)}F(\Delta_n^k(y))
  &=\partial_i\partial_j\Phi(y)
    \big(\partial_m M_j^\Delta(y_j)\partial_\ell M_i^\Delta(y_i)-\delta_{\ell,m}\partial_m M_0^\Delta(\Phi(y))\big)\\  
  &\quad  -\partial_j\Phi(y)\big(\partial_i\Phi(y)\partial_\ell\partial_m M_0^\Delta(\Phi(y))-\delta_{i,j}\partial_\ell\partial_m M_j^\Delta(y_j)\big).
}
\end{enumerate}
}

\begin{proof} (i) Let $i\in\{1,\dots,k\}$, $\ell\in\{1,\dots,n\}$, and $y\in I$ be arbitrary. Then the existence of the partial derivative $\partial_{\ell+n(i-1)}F(\Delta_n^k(y))$ and also the formula for it is a direct consequence of the standard chain rule. More precisely, 
 	\Eq{*}{
 		\partial_{\ell+n(i-1)}F(\Delta_n^k(y))
 		&=\partial_{\ell+n(i-1)}F(\Delta_n(y_1),\dots,\Delta_n(y_k))\\
 		&=\partial_i\Phi\big(M_1(\Delta_n(y_1)),\dots,M_k(\Delta_n(y_k))\big)
 		\partial_\ell M_i(\Delta_n(y_i))
 		-\partial_\ell M_0(\Phi(y),\dots,\Phi(y))\partial_i\Phi(y),
 	}
which simplifies to the formula stated in (i).

(ii) For $\alpha\in\{0,\dots,k\}$, there exists an open set $U_\alpha\subseteq I_\alpha^n$ such that $\diag(I_\alpha^n)\subseteq U_\alpha$, the first-order partial derivatives of $M_\alpha$ exist over $U_\alpha$, and their first-order partial derivatives, i.e., the second-order partial derivatives of $M_\alpha$, exist on $\diag(I_\alpha^n)$. Using the continuity of $\Phi$, by shrinking the open sets $U_1,\dots,U_k$, we can also assume
\Eq{Phi}{
  (\Phi(x_1),\dots,\Phi(x_n))\in U_0
}
provided that $x^{1}\in U_1, \dots, x^{k}\in U_k$.

Let $i,j\in\{1,\dots,k\}$, $\ell,m\in\{1,\dots,n\}$, and $y\in I$ be arbitrary. Computing the partial derivative of $F$ over $U_1\times\dots\times U_k$ with respect to its $(m+n(j-1))$th variable, i.e., with respect to the variable $x_m^j$, which is the $j$th entry of $x_m$ and the $m$th entry of $x^j$, we get
 	\Eq{*}{
 		\partial_{m+n(j-1)}F(x)
 		=\partial_j\Phi(M_1(x^1),\dots,M_k(x^k))
 		\partial_m M_j(x^j)
 		-\partial_m M_0(\Phi(x_1),\dots,\Phi(x_n))\partial_j\Phi(x_m)
 	}
for all matrices $x\in\R^{n\times k}$ with $x^1\in U_1, \dots, x^k\in U_k$. Using this equality, we can compute the partial derivative of $\partial_{m+n(j-1)}F$ at $\Delta_n^k(y)$ with respect to its $(\ell+n(i-1))$th variable, i.e., with respect to the variable $x_{\ell}^i$, which is the $i$th entry of $x_\ell$ and the $\ell$th entry of $x^i$, as follows
 	\Eq{*}{
 		\partial_{\ell+n(i-1)}&\partial_{m+n(j-1)}F(\Delta_n^k(y))\\
 		&=\partial_i\partial_j\Phi\big(M_1(\Delta_n(y_1)),\dots,M_k(\Delta_n(y_k))\big)
 		\partial_\ell M_i(\Delta_n(y_i))\partial_m M_j(\Delta_n(y_j))\\
 		&\quad+\partial_j\Phi\big(M_1(\Delta_n(y_1)),\dots,M_k(\Delta_n(y_k))\big)
 		\delta_{i,j}\partial_\ell\partial_m M_j(\Delta_n(y_j))\\  
 		&\quad-\partial_\ell\partial_m M_0(\Phi(y),\dots,\Phi(y))\partial_i\Phi(y)\partial_j\Phi(y)-\partial_m M_0(\Phi(y),\dots,\Phi(y))\delta_{\ell,m}\partial_i\partial_j\Phi(y).
 	}
Using the mean value property of $M_0,M_1,\dots,M_k$, this equality simplifies to the formula asserted in statement (ii).
\end{proof}

Our first result describes the first-order necessary condition for the validity of \eq{Mineq} in the local sense.

\Thm{PZ17}{Let $n,k\in\N$ with $n,k\geq 2$ and, for $\alpha\in\{0,\dots,k\}$, let $I_\alpha\subseteq\R$ be a nonempty open interval, $M_\alpha\colon I_\alpha^n\to I_\alpha$ be an $n$-variable mean which is partially differentiable on $\diag(I_\alpha^n)$ and let $\Phi\colon I\to I_0$ be surjective and differentiable with nonvanishing first-order partial derivatives, where $I:=I_1\times\dots\times I_k$. Assume that inequality \eq{Mineq} holds in the local sense. Then there exist constants $\lambda_1,\dots,\lambda_n\in\R_+$ such that, for all $(y_0,y)\in I_0\times I$ and $\ell\in\{1,\dots,n\}$,
\Eq{1A}{
  \lambda_\ell=\partial_\ell M_0^\Delta(y_0)
  =\partial_\ell M_1^\Delta(y_1)=\dots=\partial_\ell M_k^\Delta(y_k).
}
If, additionally, for some $\alpha\in\{0,\dots,k\}$ and $y_\alpha\in I_\alpha$, the mean $M_\alpha$ is differentiable at $\Delta_n(y_\alpha)$, then $\lambda_1+\dots+\lambda_n=1$ also holds.}

\begin{proof}
For $\alpha\in\{1,\dots,k\}$, let $U_\alpha\subseteq I_\alpha^n$ be a nonempty open set containing $\diag(I_\alpha^n)$ such that \eq{Mineq} holds for all matrices $x\in\R^{n\times k}$ with $x^1\in U_1,\dots,x^k\in U_k$. Then, according to \eq{Mineq}, $F$ is nonnegative on $U_1\times\dots\times U_k$ and, for all $y\in I$, we have $F(\Delta_n^k(y))=0$. Therefore, the first-order partial derivatives of $F$ vanish at the point $\Delta_n^k(y)$. In view of \lem{PDF}, for all $\ell\in\{1,\dots,n\}$ and $i\in\{1,\dots,k\}$, the equality $ \partial_{\ell+n(i-1)}F(\Delta_n^k(y))=0$ implies that
\Eq{*}{
  0=\partial_i\Phi(y)\big(\partial_\ell M_i^\Delta(y_i)-\partial_\ell M_0^\Delta(\Phi(y))\big).
}
Using that the partial derivatives of $\Phi$ do not vanish, for all $\ell\in\{1,\dots,n\}$, $i\in\{1,\dots,k\}$, it follows that
\Eq{ij}{
  \partial_\ell M_i^\Delta(y_i)
  =\partial_\ell M_0^\Delta(\Phi(y)).
}

We will first prove, for all $\ell\in\{1,\dots,n\}$, that the function $\partial_\ell M_0^\Delta$ is locally constant on $I_0$. Without loss of generality, we may assume that $\Phi$ is strictly increasing in its first variable. To verify the assertion, let $y_0\in I_0$ be arbitrary. Then, by the surjectivity of $\Phi$, there exists $y\in I$ such that $y_0=\Phi(y)$. Let $y_1'<y_1<y_1''$ be arbitrarily fixed elements of $I_1$. Then, by the assumed monotonicity of $\Phi$, we have 
\Eq{*}{
  y_0':=\Phi(y_1',y_2,\dots,y_k)<y_0=\Phi(y_1,y_2,\dots,y_k)<y_0'':=\Phi(y_1'',y_2,\dots,y_k).
}
Let $u\in\,]y_0',y_0''[\,$ be arbitrary. Then, by the continuity of the function $\Phi$, there exists $v\in\,]y_1',y_1''[\,$ such that $u=\Phi(v,y_2,\dots,y_k)$. Applying equality \eq{ij} for $i=2$ and $\ell\in\{1,\dots,n\}$ twice, we get
\Eq{*}{
  \partial_\ell M_0^\Delta(u)
  =\partial_\ell M_0^\Delta(\Phi(v,y_2\dots,y_k))
  =\partial_\ell M_2^\Delta(y_2)
  =\partial_\ell M_0^\Delta(\Phi(y_1,y_2\dots,y_k))
  =\partial_\ell M_0^\Delta(y_0).
}
Therefore $\partial_\ell M_0^\Delta$ is constant on $\,]y_0',y_0''[\,$, which is a neighbourhood of $y_0$. It proves that $\partial_\ell M_0^\Delta$ is differentiable at $y_0$ and $(\partial_\ell M_0^\Delta)'(y_0)=0$. The choice of $y_0$ in $I_0$ was arbitrary, hence $(\partial_\ell M_0^\Delta)'$ is identically zero on $I_0$, which is an open subinterval of $\R$. This implies that $\partial_\ell M_0^\Delta$ is constant on $I_0$. We will denote this constant by $\lambda_\ell$. Equality \eq{ij} then implies that the partial derivatives $\partial_\ell M_1^\Delta,\dots,\partial_\ell M_k^\Delta$ are also equal to the constant $\lambda_\ell$ on their domains. 

Finally, assume that, for some $\alpha\in\{0,\dots,k\}$ and $y_\alpha\in I_\alpha$, the mean $M_\alpha$ is differentiable at $(\Delta_n(y_\alpha))$. Then, by the mean value property of $M_\alpha$, we have that $M_\alpha^\Delta(y)=y$ for all $y\in I_\alpha$. Differentiating this equality with respect to $y$ at $y=y_\alpha$, we get
\Eq{*}{
  \partial_1M_i^\Delta(y_\alpha)+\dots+\partial_nM_i^\Delta(y_\alpha)=1,
}
which implies $\lambda_1+\dots+\lambda_n=1$.
\end{proof}

\Thm{PZ172}{
Let $n,k\in\N$ with $n,k\geq 2$ and, for $\alpha\in\{0,\dots,k\}$, let $I_\alpha\subseteq\R$ be a nonempty open interval, $M_\alpha\colon I_\alpha^n\to I$ be an $n$-variable mean which is twice differentiable on $\diag(I_\alpha^n)$ and let $\Phi\colon I\to I_0$ be surjective and twice differentiable with nonvanishing first-order partial derivatives, where $I:=I_1\times\dots\times I_k$. Assume that inequality \eq{Mineq} holds in the local sense. Then there exist constants $\lambda_1,\dots,\lambda_n\in\R_+$ with $\lambda_1+\dots+\lambda_n=1$ such that, for all $(y_0,y)\in I_0\times I$ and $\ell\in\{1,\dots,n\}$, the equalities in \eq{1A} hold. In addition, for all $y\in I$, the $(nk)\times(nk)$ matrix whose $(\ell+n(i-1),m+n(j-1))$th entry, where $i,j\in\{1,\dots,k\}$ and $\ell,m\in\{1,\dots,n\}$, is given by 
\Eq{Aij}{
	\partial_i\partial_j\Phi(y)
	(\lambda_m\lambda_\ell-\delta_{\ell,m}\lambda_m)
	-\partial_i\Phi(y)\partial_j\Phi(y)\partial_\ell\partial_m M_0^\Delta(\Phi(y))+\delta_{i,j}\partial_j\Phi(y)\partial_\ell\partial_m M_j^\Delta(y_j)
    }
is positive semidefinite.
}

\begin{proof}
For $\alpha\in\{1,\dots,k\}$, let $U_\alpha\subseteq I_\alpha^n$ be a nonempty open set containing $\diag(I_\alpha^n)$ such that \eq{Mineq} and \eq{Phi} hold for all matrices $x\in\R^{n\times k}$ with $x^1\in U_1,\dots,x^k\in U_k$. Then, using \thm{PZ17}, condition \eq{1A} is valid with some nonnegative constants $\lambda_1,\dots,\lambda_n$ satisfying also $\lambda_1+\dots+\lambda_n=1$.
	
According to \eq{Mineq}, $F$ is nonnegative on $U_1\times\dots\times U_k$ and, for all $y\in I$, we have $F(\Delta_n^k(y))=0$, that is, $F$ has a (local) minimum at $\Delta_n^k(y)$. Therefore, its second derivative, i.e., the $(nk)\times(nk)$ matrix $\big(\partial_\alpha\partial_\beta F(\Delta_n^k(y))\big)$ is positive semidefinite. In view of \lem{PDF}, for all $y\in I_1\times\dots\times I_k$, it follows that the $(nk)\times(nk)$ matrix whose $(\ell+n(i-1),m+n(j-1))$th entry, where $i,j\in\{1,\dots,k\}$ and $\ell,m\in\{1,\dots,n\}$, is given by 
\Eq{*}{
	\partial_i\partial_j\Phi(y)\big(\partial_m M_j^\Delta(y_j)\partial_\ell M_i^\Delta(y_i)-\delta_{\ell,m}\partial_m M_0^\Delta(\Phi(y))\big)
		-\partial_j\Phi(y)\big(\partial_i\Phi(y)\partial_\ell\partial_m M_0^\Delta(\Phi(y))-\delta_{i,j}\partial_\ell\partial_m M_j^\Delta(y_j)\big)
}
is positive semidefinite. Applying the equalities from  \eq{1A}, the statement follows.
\end{proof}

\section{Hölder- and Minkowski-type inequalities for nonsymmetric generalized Bajraktarević means}

In order to apply the results from the previous section for nonsymmetric generalized Bajraktarević means, we need to compute their partial derivatives on $\diag(I^n)$. For this aim, we recall the following result, which was obtained by the authors in \cite{GruPal20}. 

\Lem{DB}{
Let $n,k\in\N, d\in\{1,2\}$, let $f\colon  I\to\R$ be a $d$ times differentiable function with a nonvanishing first derivative, $p=(p_1,\dots,p_n)\colon I\to\R_{+}^n$ and set $p_0:=p_1+\dots+p_n$. Then we have the following assertions.
\begin{enumerate}[(i)]
	\item If $d=1$ and $p$ is continuous, then, for all $\ell\in \{1,\dots,n\}$, the first-order partial derivative $\partial_\ell A_{f,p}$ exists on $\diag(I^n)$ and
	\Eq{*}{
		\partial_\ell A_{f,p}^\Delta=\frac{p_\ell}{p_0}.
		}
	\item If $d=2$ and $p$ is continuously differentiable, then, for all $\ell,m\in \{1,\dots,n\}$, the second-order partial derivatives $\partial_\ell^2 A_{f,p}$ and $\partial_\ell\partial_m A_{f,p}$ exist on $\diag(I^n)$ and
	\Eq{*}{
		\partial_\ell^2A_{f,p}^\Delta=2\frac{p_\ell'(p_0-p_\ell)}{p_0^2}+\frac{p_\ell(p_0-p_\ell)}{p_0^2}\cdot\frac{f''}{f'},\qquad
		\partial_\ell\partial_m A_{f,p}^\Delta =-\frac{(p_\ell p_m)'}{p_0^2}-\frac{p_\ell p_m}{p_0^2}\cdot\frac{f''}{f'}\qquad(\ell\neq m).
	}
\end{enumerate}
}

\Thm{PZ172b}{
Let $n,k\in\N$ with $n,k\geq 2$ and, for $\alpha\in\{0,\dots,k\}$, let $I_\alpha\subseteq\R$ be a nonempty open interval, $f_\alpha\colon I_\alpha\to\R$ be a  differentiable function with a nonvanishing first derivative and let $p^\alpha=(p_1^\alpha,\dots,p_n^\alpha) \colon I_\alpha\to\R_+^n$ be continuous, set $p_0^\alpha:=p_1^\alpha+\dots+p_n^\alpha$ and denote $I:=I_1\times\dots\times I_k$. Let $\Phi\colon I\to I_0$ be surjective and differentiable with nonvanishing first-order partial derivatives. Assume that inequality \eq{ineq} holds in the local sense. Then there exist constants $\lambda_1,\dots,\lambda_n\in\R_+$ with $\lambda_1+\dots+\lambda_n=1$ such that, for all $\alpha\in\{0,\dots,k\}$ and $\ell\in\{1,\dots,n\}$, 
	\Eq{1C}{
	  p_\ell^\alpha=\lambda_\ell p_0^\alpha
	}
holds on $I_\alpha$. If, additionally, for $\alpha\in\{0,\dots,k\}$, $f_\alpha$ is twice differentiable, $p^\alpha$ is continuously differentiable and $\Phi$ is twice differentiable, then the $k\times k$ matrix $\Gamma(y)$ given by 
    \Eq{1D}{
    \Gamma_{i,j}(y):=\bigg(&-\partial_i\partial_j\Phi(y)
    -\partial_j\Phi(y)\partial_i\Phi(y)\bigg(2\frac{(p^0_0)'}{p^0_0}+\frac{f_0''}{f_0'}\bigg)(\Phi(y))
    +\delta_{i,j}\partial_j\Phi(y)\bigg(2\frac{(p^j_0)'}{p^j_0}+\frac{f_j''}{f_j'}\bigg)(y_j)\bigg)_{i,j=1}^k
    }
is positive semidefinite for all $y\in I$.
}

\begin{proof} For $\alpha\in\{0,\dots,k\}$, let $M_\alpha=A_{f_\alpha,p^\alpha}$ and apply \thm{PZ172} to this setting. Then $M_\alpha$ is partially differentiable on $\diag(I_\alpha^n)$ and inequality \eq{Mineq} holds in the local sense. According to the first assertion of \thm{PZ172} and by the first statement of \lem{DB}, there exist $\lambda_1,\dots,\lambda_n\in\R_+$ such that, for all $(y_0,y)\in I_0\times I$ and $\ell\in\{1,\dots,n\}$, the equalities in \eq{1A} hold, i.e., for $\alpha\in\{0,\dots,k\}$,
\Eq{*}{
  \frac{p^\alpha_\ell}{p^\alpha_0}(y_\alpha)
  =\partial_\ell A_{f_\alpha,p^\alpha}^\Delta(y_\alpha)=\lambda_\ell.
}
This shows that \eq{1C} is valid on $I_\alpha$ for all $\alpha\in\{0,\dots,k\}$ and $\ell\in\{1,\dots,n\}$.
In view of the definition of $p_0^\alpha$, these equalities imply that $\lambda_1+\dots+\lambda_n=1$ is also valid.

Assume now that, additionally, for $\alpha\in\{0,\dots,k\}$, $f_\alpha$ is twice differentiable, $p^\alpha$ is continuously differentiable and $\Phi$ is twice differentiable.
Using \eq{1C}, according to the second assertion of \lem{DB}, we have that
\Eq{*}{
	\partial_\ell^2A_{f_\alpha,p^\alpha}^\Delta &=2\frac{(p^\alpha_\ell)'(p^\alpha_0-p^\alpha_\ell)}{(p^\alpha_0)^2}+\frac{p^\alpha_\ell(p^\alpha_0-p^\alpha_\ell)}{(p^\alpha_0)^2}\cdot\frac{f''_\alpha}{f'_\alpha}=\lambda_\ell(1-\lambda_\ell)\bigg(2\frac{(p^\alpha_0)'}{p^\alpha_0}+\frac{f''_\alpha}{f'_\alpha}\bigg),\\
	\partial_\ell\partial_m A_{f_\alpha,p^\alpha}^\Delta &=-\frac{(p^\alpha_\ell p^\alpha_m)'}{(p^\alpha_0)^2}-\frac{p^\alpha_\ell p^\alpha_m}{(p^\alpha_0)^2}\cdot\frac{f''_\alpha}{f'_\alpha}=-\lambda_\ell\lambda_m\bigg(2\frac{(p^\alpha_0)'}{p^\alpha_0}+\frac{f''_\alpha}{f'_\alpha}\bigg)\qquad(\ell\neq m).
}
Therefore, for all $\alpha\in\{0,\dots,k\}$ and $\ell,m\in\{1,\dots,n\}$,
\Eq{lm}{
  \partial_\ell\partial_m A_{f_\alpha,p^\alpha}^\Delta
  =\lambda_m(\delta_{\ell,m}-\lambda_\ell)\bigg(2\frac{(p^\alpha_0)'}{p^\alpha_0}+\frac{f''_\alpha}{f'_\alpha}\bigg).
}

By the second assertion of \thm{PZ172}, for all $y\in I$, the $(nk)\times (nk)$ matrix whose $(\ell+n(i-1),m+n(j-1))$th entry, where $i,j\in\{1,\dots,k\}$ and $\ell,m\in\{1,\dots,n\}$, is given by \eq{Aij} is positive semidefinite. Using formula \eq{lm}, we can conclude that the matrix whose $(\ell+n(i-1),m+n(j-1))$th entry is given by
\Eq{*}{
 \lambda_m(\delta_{\ell,m}-\lambda_\ell)
 \Gamma_{i,j}(y)
}
is positive semidefinite.

If a matrix is positive semidefinite, then every minor of the matrix is also positive semidefinite. Therefore, the $k\times k$ submatrix with entries $(1+n(i-1),1+n(j-1))$th, where $i,j\in\{1,\dots,k\}$, is also positive semidefinite, which implies the statement.
\end{proof}

In the next result, we reformulate the positive semidefiniteness condition from the above theorem in terms of a convexity property.

\Thm{PZ173}{
Let $n,k\in\N$ with $n,k\geq 2$ and, for $\alpha\in\{0,\dots,k\}$, let $I_\alpha\subseteq\R$ be a nonempty open interval, $f_\alpha\colon I_\alpha\to\R$ be a twice differentiable function with a nonvanishing first derivative and let $p^\alpha=(p_1^\alpha,\dots,p_n^\alpha) \colon I_\alpha\to\R_+^n$ be continuously differentiable, set $p_0^\alpha:=p_1^\alpha+\dots+p_n^\alpha$ and denote $I:=I_1\times\dots\times I_k$. Let $\Phi\colon I\to I_0$ be surjective and twice differentiable with nonvanishing first-order partial derivatives. Assume that inequality \eq{ineq} holds in the local sense. 
Finally, for $\alpha\in\{0,\dots,k\}$, define the function $\varphi_\alpha\colon I_\alpha\to\R$ and then $\varphi\colon I\to\R^k$ by
\Eq{*}{
  \varphi_\alpha:=\int (p_0^\alpha)^2f_\alpha'
  \qquad\mbox{and}\qquad
  \varphi(y):=(\varphi_1(y_1),\dots,\varphi_k(y_k)).
}
Then $\varphi_1,\dots,\varphi_k$ and $\varphi$ are twice differentiable and invertible functions and the map $\Psi\colon \varphi(I)\to\R$ defined by
\Eq{*}{
  \Psi(u):=\varphi_0(\Phi(\varphi^{-1}(u)))
}
is concave if $f_0'>0$ and convex if $f_0'<0$.}

\begin{proof} According to \thm{PZ172b}, our assumptions imply that the matrix-valued map $\Gamma\colon I\to\R^{k\times k}$ defined by \eq{1D} has positive semidefinite values.

Without loss of generality, we can assume that $f_0'>0$. Let $\alpha\in\{0,\dots,k\}$. Then the integrand in the definition of $\varphi_\alpha$ is either positive everywhere or negative everywhere, therefore $\varphi_\alpha$ is a twice differentiable function with a nonvanishing first derivative, hence it is strictly monotone and it has a twice differentiable inverse $\varphi_\alpha^{-1}\colon \varphi_\alpha(I_\alpha)\to I_\alpha$. Furthermore, we have that
\Eq{fij}{
  \frac{\varphi_\alpha''}{\varphi_\alpha'}
  =\frac{((p_0^j)^2f_\alpha')'}{(p_0^j)^2f_\alpha'}
  =\frac{2p_0^j(p_0^j)'f_\alpha'+(p_0^j)^2f_\alpha''}{(p_0^j)^2f_\alpha'}
  =2\frac{(p^j_0)'}{p^j_0}+\frac{f_\alpha''}{f_\alpha'}
  \qquad(\alpha\in\{0,\dots,k\}).
}
It follows from the definition of $\varphi$ that
\Eq{*}{
  \varphi^{-1}(u)=(\varphi_1^{-1}(u_1),\dots,\varphi_k^{-1}(u_k))
  \qquad(u\in\varphi_1(I_1)\times\dots\times \varphi_k(I_k)).
}
Thus, it is clear that $\varphi$ and its inverse are
also twice differentiable maps. 

In order to show that $\Psi$ is concave, we will prove that $\Psi''$ is negative semidefinite over $\varphi(I)$. First, we compute the first and then the second-order partial derivatives of $\Psi$. For $i,j\in\{0,\dots,k\}$ and $u\in\varphi(I)$, using standard calculus rules, we obtain
\Eq{*}{
  \partial_j \Psi(u)
  &=\varphi_0'(\Phi(\varphi^{-1}(u)))\cdot\partial_j\Phi(\varphi^{-1}(u))\cdot\frac{1}{\varphi_j'(\varphi_j^{-1}(u_j))}
}
and
\Eq{*}{
  \partial_i\partial_j \Psi(u)
  &=\varphi_0''(\Phi(\varphi^{-1}(u)))\cdot\partial_i\Phi(\varphi^{-1}(u))\cdot\partial_j\Phi(\varphi^{-1}(u))
  \cdot\frac{1}{\varphi_i'(\varphi_i^{-1}(u_i))}
  \cdot\frac{1}{\varphi_j'(\varphi_j^{-1}(u_j))}\\
  &\quad+\varphi_0'(\Phi(\varphi^{-1}(u)))\cdot\partial_i\partial_j\Phi(\varphi^{-1}(u))\cdot\frac{1}{\varphi_i'(\varphi_i^{-1}(u_i))}\cdot\frac{1}{\varphi_j'(\varphi_j^{-1}(u_j))}\\
  &\quad-\delta_{i,j}\varphi_0'(\Phi(\varphi^{-1}(u)))\cdot\partial_j\Phi(\varphi^{-1}(u))\cdot\frac{\varphi_j''(\varphi_j^{-1}(u_j))}{\varphi_j'(\varphi_j^{-1}(u_j))^3}\\
  &=\frac{\varphi_0'(\Phi(\varphi^{-1}(u)))}{\varphi_i'(\varphi_i^{-1}(u_i))\cdot\varphi_j'(\varphi_j^{-1}(u_j))}\big(\partial_i\Phi(\varphi^{-1}(u))\cdot\partial_j\Phi(\varphi^{-1}(u))\cdot
  \frac{\varphi_0''}{\varphi_0'}(\Phi(\varphi^{-1}(u))\\
  &\hspace{5.2cm}+\partial_i\partial_j\Phi(\varphi^{-1}(u))
  -\delta_{i,j}\partial_j\Phi(\varphi^{-1}(u))\cdot\frac{\varphi_j''}{\varphi_j'}(\varphi_j^{-1}(u_j))\big).
}
Now using the equalities in \eq{fij} and \eq{1D}, it follows that
\Eq{*}{
  \partial_i\partial_j \Psi(u)
  =\frac{\varphi_0'(\Phi(\varphi^{-1}(u)))}{\varphi_i'(\varphi_i^{-1}(u_i))\cdot\varphi_j'(\varphi_j^{-1}(u_j))}\cdot(-\Gamma_{i,j}(\varphi^{-1}(u))).
}

Therefore, for all $u\in\varphi(I)$, we obtain that $\Psi''(u)=(\partial_i\partial_j \Psi(u))_{i,j=1}^k$ is negative semidefinite. This implies that $\Psi$ is concave on $\varphi(I)$.
\end{proof}

\Rem{17PZ}{It can be seen from the above argument that the concavity of the auxiliary function $\Psi$ is not merely a consequence of the positive semidefiniteness of the matrix-valued function $\Gamma$ but, in fact, it is equivalent to it. On the other hand, if all the weight functions are equal to constant 1, then $\varphi_\alpha=f_\alpha$ and, in this case, according to the theory of quasiarithmetic means (cf.\ \cite{HarLitPol34}), the concavity of the function $\Psi$ is also sufficient for inequality \eq{ineq} to be valid in the global sense.}

The following results establish sufficient conditions for inequality \eq{ineq} to be valid in the local as well as in the global sense.

\Thm{17PZloc}{Let $k\in\N$ and, for $\alpha\in\{0,\dots,k\}$, let $I_\alpha\subseteq\R$ be a nonempty open interval, $f_\alpha\colon I_\alpha\to\R$ be differentiable with a nonvanishing derivative, $p_0^\alpha\colon I_\alpha\to\R_+$ and denote $I:=I_1\times\dots\times I_k$. Furthermore, let $\Phi\colon I\to I_0$ be partially differentiable. Assume that there exists an open set $V\subseteq I^2$ containing $\diag(I^2)$ such that, for all $(u,y)\in V$, the inequality 
\Eq{GSC0}{
  \frac{p_0^0(\Phi(y))(f_0(\Phi(y))-f_0(\Phi(u)))}{p_0^0(\Phi(u))f'_0(\Phi(u))}
  \leq \sum_{j=1}^k\partial_j\Phi(u)
  \frac{p_0^j(y_j)(f_j(y_j)-f_j(u_j))}{p_0^j(u_j)f_j'(u_j)}
}
holds. Then, for all $n\in\N$ and $\lambda\in\R_+^n$ with $\lambda_1+\dots+\lambda_n=1$, the inequality 
\Eq{MineqBL}{
	A_{f_0,p_0^0\lambda}(\Phi(x_1),\dots,\Phi(x_n)) &\leq \Phi(A_{f_1,p_0^1\lambda}(x^1),\dots,A_{f_k,p_0^k\lambda}(x^k))
}
is valid in the local sense.}

\begin{proof} Let $n\in\N$, $\lambda\in\R_+^n$ with $\lambda_1+\dots+\lambda_n=1$ be fixed and construct the set $U\subseteq I^n$ as follows:
\Eq{U}{
  U:=\bigcap_{i=1}^n\{x\in\R^{n\times k}\colon \transp{x}\in I^n, \,\, (A_{f_1,p_0^1\lambda}(x^1),\dots,A_{f_k,p_0^k\lambda}(x^k),x_i^1,\dots,x_i^k)\in V\}.
}
Then, due to the continuity of the mean $A_{f_\alpha,p_0^\alpha\lambda}$, each member of the intersection is open, and hence so is $U$. On the other hand, if $\transp{x}\in\diag(I^n)$, then, for all $\alpha\in\{1,\dots,k\}$, we have $x_1^\alpha=\dots=x_n^\alpha=A_{f_\alpha,p_0^\alpha\lambda}(x^\alpha)$, whence, by the properties of $V$, $(A_{f_1,p_0^1\lambda}(x^1),\dots,A_{f_k,p_0^k\lambda}(x^k),x_i^1,\dots,x_i^k)\in V$ holds for all $i\in\{1,\dots,n\}$. This shows that $U$ contains $\diag(I^n)$. 

We now prove that, for all $x\in\R^{n\times k}$ with $\transp{x}\in U$, inequality \eq{MineqBL} is valid. Let us define, for $\alpha\in\{1,\dots,k\}$,
\Eq{*}{
  y_\alpha:=A_{f_\alpha,p_0^\alpha\lambda}(x^\alpha)
}
and set $u\in I$. As a consequence of this definition, it follows that
\Eq{eqs}{
  \sum_{i=1}^n
  \lambda_i p_0^\alpha(x_i^\alpha)(f_\alpha(x_i^\alpha)-f_\alpha(u_\alpha))=0
  \qquad(\alpha\in\{1,\dots,k\}).
}
On the other hand, for all $\alpha\in\{1,\dots,k\}$ and $i\in\{1,\dots,n\}$, we have that $(u_1,\dots,u_k,x_i^1,\dots,x_i^k)\in V$. Therefore, we can apply \eq{GSC0} with $(u_1,\dots,u_k,y_1,\dots,y_k) :=(u_1,\dots,u_k,x_i^1,\dots,x_i^k)$. Then multiplying each inequality by $\lambda_i$, summing up the inequalities so obtained, and using the identities in \eq{eqs}, we obtain
\Eq{*}{
  \sum_{i=1}^n\frac{\lambda_i p_0^0(\Phi(x_i))(f_0(\Phi(x_i)-f_0(\Phi(u)))}{p_0^0(\Phi(u))f'_0(\Phi(u))} &\leq \sum_{i=1}^n\sum_{\alpha=1}^k\partial_\alpha\Phi(u)
  \frac{\lambda_i p_0^\alpha(x_i^\alpha)(f_\alpha(x_i^\alpha)-f_\alpha(u_\alpha))}{p_0^\alpha(u_\alpha)f_\alpha'(u_\alpha)}\\& =\sum_{\alpha=1}^k\frac{\partial_\alpha\Phi(u)}{p_0^\alpha(u_\alpha)f_\alpha'(u_\alpha)}\sum_{i=1}^n
  \lambda_i p_0^\alpha(x_i^\alpha)(f_\alpha(x_i^\alpha)-f_\alpha(u_\alpha))=0.
}
Therefore,
\Eq{*}{
  \sum_{i=1}^n&\frac{\lambda_i p_0^0(\Phi(x_i))(f_0(\Phi(x_i)-f_0(\Phi(u)))}{p_0^0(\Phi(u))f'_0(\Phi(u))}\leq 0.
}
Assume that $f_0'$ is positive. Then $f_0$ is strictly increasing and the above inequality is equivalent to
\Eq{I1}{
  \sum_{i=1}^n&\lambda_i p_0^0(\Phi(x_i))(f_0(\Phi(x_i)-f_0(\Phi(u)))\leq 0.
}
Rearranging this inequality, we obtain
\Eq{I2}{
  \frac{\sum_{i=1}^n\lambda_i p_0^0(\Phi(x_i))f_0(\Phi(x_i))}{\sum_{i=1}^n\lambda_i p_0^0(\Phi(x_i))}&\leq f_0(\Phi(u)).
}
Applying $f_0^{-1}$ side by side and using that $f_0^{-1}$ is strictly increasing, we can conclude that
\Eq{*}{
  A_{f_0,p_0^0\lambda}(\Phi(x_1),\dots,\Phi(x_n))
  &=f_0^{-1}\bigg(\frac{\sum_{i=1}^n\lambda_i p_0^0(\Phi(x_i))f_0(\Phi(x_i))}{\sum_{i=1}^n\lambda_i p_0^0(\Phi(x_i))}\bigg)\\
  &\leq \Phi(u)=\Phi(A_{f_1,p_0^1\lambda}(x^1),\dots,A_{f_k,p_0^k\lambda}(x^k)),
}
which completes the proof of inequality \eq{MineqBL}. In the case, when $f'_0$ is everywhere negative, the inequalities \eq{I1} and \eq{I2} are reversed, however
$f_0^{-1}$ is strictly decreasing, thus we arrive at the same conclusion.
\end{proof}

\Rem{17PZloc}{
	In view of \thm{PZ172b}, the weight functions of the generalized nonsymmetric Bajraktarevi\'c means appearing in \eq{ineq} necessarily are of the form given by \eq{1C}. Therefore, the local as well as the global validity of \eq{ineq} immediately follows from the local as well as the global validity of \eq{MineqBL}, respectively.	
}

\Thm{17PZgl}{Let $k\in\N$ and, for $\alpha\in\{0,\dots,k\}$, let $I_\alpha\subseteq\R$ be a nonempty open interval, $f_\alpha\colon I_\alpha\to\R$ be differentiable with a nonvanishing derivative, $p_0^\alpha\colon I_\alpha\to\R_+$ and let $\Phi\colon I\to I_0$ be partially differentiable, where $I:=I_1\times\dots\times I_k$. Assume, for all $u,y\in I$, that inequality \eq{GSC0} is satisfied. Then, for all $n\in\N$, $\lambda\in\R_+^n$ with $\lambda_1+\dots+\lambda_n=1$, inequality \eq{MineqBL} holds in the global sense.}

\begin{proof} If \eq{GSC0} is satisfied for all $u,y\in I$, then the condition of the previous theorem is validated with the open set $V:=I^2$ and the open set $U$ constructed by \eq{U} equals $I^n$. Hence inequality \eq{MineqBL} holds for all $x\in\R^{n\times k}$ with $\transp{x}\in U$, i.e., it holds in the global sense.
\end{proof}

The next result establishes a necessary condition for \eq{GSC0} to be satisfied in the local sense.

\Thm{NecC}{
	Let $k\in\N$ and, for $\alpha\in\{0,\dots,k\}$, let $I_\alpha\subseteq\R$ be a nonempty open interval, $f_\alpha\colon I_\alpha\to\R$ be twice differentiable with a nonvanishing first derivative and $p_0^\alpha\colon I_\alpha\to\R_+$ be twice differentiable. In addition, let $\Phi\colon I\to I_0$ be twice differentiable, where $I:=I_1\times\dots\times I_k$. Assume that there exists an open set $V\subseteq I^2$ with $\diag(I^2)\subseteq V$ such that \eq{GSC0} is satisfied for all $(u,y)\in V$. Then the matrix-valued function $\Gamma\colon I\to\R^{k\times k}$ defined by \eq{1D} takes positive semidefinite values.
}

\begin{proof} If \eq{GSC0} is satisfied for all $(u,y)\in V$ then, for all fixed $y\in I$, the map $\Psi_y\colon I\to\R$ defined as
\Eq{Psi_y}{
  \Psi_y(u):=\sum_{\alpha=1}^k\partial_\alpha\Phi(y)
  \frac{p_0^\alpha(u_\alpha)(f_\alpha(u_\alpha)-f_\alpha(y_\alpha))}{p_0^\alpha(y_\alpha)f_\alpha'(y_\alpha)}-\frac{p_0^0(\Phi(u))(f_0(\Phi(u))-f_0(\Phi(y)))}{p_0^0(\Phi(y))f'_0(\Phi(y))}
}
has a local minimum at $u=y$. This function is twice differentiable according to our assumptions. Therefore, the second derivative matrix of it at $u=y$ is positive semidefinite. For $i\in\{1,\dots,k\}$, we have
\Eq{DiPsi}{
    \partial_i\Psi_y(u)
    &=\frac{\partial_i\Phi(y)}{p_0^i(y_i)f_i'(y_i)}\big((p_0^i)'(u_i)(f_i(u_i)-f_i(y_i))+p_0^i(u_i)f_i'(u_i)\big) \\
    &\quad-\frac{\partial_i\Phi(u)}{p_0^0(\Phi(y))f_0'(\Phi(y))}\big((p_0^0)'(\Phi(u))(f_0(\Phi(u))-f_0(\Phi(y)))+p_0^0(\Phi(u))f_0'(\Phi(u))\big).
}
Thus, for $i,j\in\{1,\dots,k\}$, we obtain
\Eq{DijPsi}{
    \partial_i\partial_j\Psi_y(u)
    &=\delta_{ij}\frac{\partial_j\Phi(y)}{p_0^j(y_j)f_j'(y_j)}\big((p_0^j)''(u_j)(f_j(u_j)-f_j(y_j))+2(p_0^j)'(u_j)f_j'(u_j)+p_0^j(u_j)f_j''(u_j)\big) \\
    &\quad - 
    \frac{\partial_i\partial_j\Phi(u)}{p_0^0(\Phi(y))f_0'(\Phi(y))}\big(
    (p_0^0)'(\Phi(u))(f_0(\Phi(u))-f_0(\Phi(y)))+ p_0^0(\Phi(u))f_0'(\Phi(u))\big)\\
    &\quad - 
    \frac{\partial_i\Phi(u)\partial_j\Phi(u)}{p_0^0(\Phi(y))f_0'(\Phi(y))}\big((p_0^0)''(\Phi(u))(f_0(\Phi(u))-f_0(\Phi(y))) \\&
    \hspace{4cm}+2(p_0^0)'(\Phi(u))f_0'(\Phi(u)) 
    +p_0^0(\Phi(u))f_0''(\Phi(u))\big).
}
Hence, after substituting $u:=y$ in the above equality, we get
\Eq{PsiGam}{
    \partial_i\partial_j\Psi_y(y)=\delta_{ij}\partial_j\Phi(y)\bigg( 2\frac{(p_0^j)'}{p_0^j}+\frac{f_j''}{f_j'}\bigg)(y_j)-\partial_i\partial_j\Phi(y)
    -\partial_i\Phi(y)\partial_j\Phi(y)\bigg(2\frac{(p_0^0)'}{p_0^0}+\frac{f_0''}{f_0'}\bigg)(\Phi(y))=\Gamma_{i,j}(y).
}
This shows the pointwise positive semidefiniteness of the matrix-valued function $\Gamma\colon I\to\R^{k\times k}$.
\end{proof}

\Thm{SP}{Let $k\in\N$ and, for $\alpha\in\{0,\dots,k\}$, let $I_\alpha\subseteq\R$ be a nonempty open interval, $f_\alpha\colon I_\alpha\to\R$ be twice continuously differentiable with a nonvanishing first derivative, $p_0^\alpha\colon I_\alpha\to\R_+$ be twice continuously differentiable and denote $I:=I_1\times\dots\times I_k$. Furthermore, let $\Phi\colon I\to I_0$ be twice continuously differentiable. Assume, for all $y\in I$, that the $k\times k$ matrix $\Gamma(y)$ is positive definite. Then, for all $n\in\N$, $\lambda\in\R_+^n$ with $\lambda_1+\dots+\lambda_n=1$, inequality \eq{MineqBL} holds in the local sense.}

\begin{proof} For all fixed $y\in I$, define the function $\Psi_y\colon I\to\R$ by the formula \eq{Psi_y}. This function is twice continuously differentiable according to our assumptions and $\Psi_y(y)=0$ holds for all $y\in I$. After simple computations, for $i,j\in\{1,\dots,k\}$ and $u\in I$, we obtain that the equalities \eq{DiPsi} and \eq{DijPsi} hold.

Putting $u:=y$ into equality \eq{DiPsi}, we can see that
\Eq{*}{
  \partial_i\Psi_y(y)=0 \qquad(i\in\{1,\dots,k\}),
}
that is, $\Psi'_y(y)=\big(\partial_i\Psi_y(y)\big)_{i=1}^k=0$ holds for all $y\in I$.

Substituting $u:=y$ into equality \eq{DijPsi}, we can conclude that \eq{PsiGam} is valid for all $y\in I$.
According to the positive definiteness of the matrix-valued function $\Gamma\colon I\to\R^{k\times k}$, it follows that $\Psi_y''(y) :=\big(\partial_i\partial_j\Psi_y(y)\big)_{i,j=1}^k$ is positive definite for all $y\in I$.

For $u,y\in I$, denote the smallest eigenvalue of the $k\times k$ symmetric matrix $\Psi_y''(u) :=\big(\partial_i\partial_j\Psi_y(u)\big)_{i,j=1}^k$ by $\psi(u,y)$. In view of our twice continuous differentiability assumptions, the map $(u,y)\mapsto \Psi_y''(u)$ is continuous, therefore the map $(u,y)\mapsto\psi(u,y)$ is also continuous. On the other hand, for all $y\in I$, we have that $\Psi_y''(y)=\Gamma(y)$ is positive definite, which implies that $\psi(y,y)>0$. Therefore, there exists an open set $W\subseteq I^2$ containing $\diag(I^2)$ on which $\psi$ is positive. 

By the Taylor Mean Value Theorem, for all $u,y\in I$, there exists $t\in[0,1]$ such that
\Eq{Psi}{
  \Psi_y(u)&=\Psi_y(y)+\Psi_y'(y)(y-u)+\frac{1}{2}(y-u)^T\Psi_y''(tu+(1-t)y)(y-u)\\
  &=\frac{1}{2}(y-u)^T\Psi_y''(tu+(1-t)y)(y-u)
  \geq \frac{1}{2}\psi(tu+(1-t)y,y)\|y-u\|^2.
}

Define
\Eq{*}{
  V:=\{(u,y)\in I^2\mid [y,u]\times\{y\}=[(y,y),(u,y)]\subseteq W\}.
}
We will show that $V$ is an open subset of $I^2$. To see this, let $(u,y)\in V$ be arbitrary. Since the segment $[(y,y),(u,y)]$ is a compact subset of $W$, it follows that it is disjoint from $W^c:=\R^{2k}\setminus W$, which is the complement of $W$ and hence it is a closed set. Therefore, there exists a positive number $r$ so that the distance of every point of the segment $[(y,y),(u,y)]$ from $W^c$ is at least $r$. Let $v,x\in\R^k$ such that $\|v-u\|<r/2$ and $\|x-y\|<r/2$ hold. Let $(w,x)$ be an arbitrary point of the segment $[(x,x),(v,x)]$. Then there exists $t\in[0,1]$ such that $w=tv+(1-t)x$.
Therefore,
\Eq{*}{
  \|w-(tu+(1-t)y)\|\leq t\|v-u\|+(1-t)\|x-y\|<r/2,
}
which implies that
\Eq{*}{
 \|(w,x)-(tu+(1-t)y,y)\|
 &\leq \|(w,x)-(tu+(1-t)y,x)\|+\|(tu+(1-t)y,x)-(tu+(1-t)y,y)\|\\
 &=\|w-(tu+(1-t)y)\|+\|x-y\|<r.
}
In other words, any point of the segment $[(x,x),(v,x)]$ is closer to some point of the segment $[(y,y),(u,y)]$ than $r$. This yields that the segment $[(x,x),(v,x)]\subseteq W$, i.e., $(v,x)\in V$ whenever $\|v-u\|<r/2$ and $\|x-y\|<r/2$ hold, consequently $(u,y)$ is an interior point of $V$. This completes the proof of the openness of $V$. On the other hand, it is obvious that $V$ also contains $\diag(I^2)$.

Using that $\psi$ is positive on $W$, it follows from \eq{Psi} that $\Psi_y(u)\geq0$ for all $(u,y)\in V$. Therefore, by applying \thm{17PZloc}, it follows that \eq{MineqBL} holds in the local sense.
\end{proof}

\section{Inequalities for Gini means}

In this section, we apply the above results to some important particular cases of \eq{ineq}. First, we deal with cases obtained by specializing the means $M_\alpha$ in \eq{ineq} for all $\alpha\in\{0,\dots,k\}$. Then we draw some conclusions by choosing the function $\Phi$ in \eq{ineq} to be the map $k$-variable addition and multiplication. 

Let us recall the definition of the \emph{weighted $n$-variable Hölder (or power) mean of parameter $r\in\R$ and weight vector $\lambda\in\R_+^n$ with $\lambda_1+\dots+\lambda_n=1$}, and the \emph{weighted $n$-variable Gini mean corresponding to the pair parameters $(r,s)\in\R^2$ and weight vector $\lambda\in\R_+^n$ with $\lambda_1+\dots+\lambda_n=1$}:
\Eq{*}{
	H_{r;\lambda}(x_1,\dots,x_n)&:=\begin{cases}
		\big(\lambda_1x_1^r+\dots+\lambda_nx_n^r\big)^\frac{1}{r} & \mbox{if }r\neq 0, \\[1mm]
		x_1^{\lambda_1}\cdots x_n^{\lambda_n} & \mbox{if } r=0;
	\end{cases} \\[2mm]
	G_{r,s;\lambda}(x_1,\dots,x_n)&:=\begin{cases}
		\bigg(\dfrac{\lambda_1x_1^r+\dots+\lambda_nx_n^r}{\lambda_1x_1^s+\dots+\lambda_nx_n^s}\bigg)^\frac{1}{r-s} & \mbox{if }r\neq s, \\[4mm]
		\exp\bigg(\dfrac{\lambda_1x_1^r\ln(x_1)+\dots+\lambda_nx_n^r\ln(x_n)}{\lambda_1x_1^r+\dots+\lambda_nx_n^r}\bigg) & \mbox{if } r=s.
	\end{cases}
}
It is clear that in the particular case $q=0$, the mean $G_{p,q;\lambda}$ simplifies to $H_{p;\lambda}$.

For $(r,s)\in\R^2$ we also define the function $\chi_{r,s}\colon \R_+\to\R$ by
\Eq{*}{
  \chi_{r,s}(t)
  :=\begin{cases}
   \dfrac{t^r-t^s}{r-s} & \mbox{if } r\neq s, \\[3mm]
   t^r\ln(t) & \mbox{if } r=s.
   \end{cases}
}
Then, for all $\ell\in\{1,\dots,n\}$, with $p_\ell(t):=\lambda_\ell t^s$ and $f(t):=t^{r-s}$ if $r\neq s$ or $f(t):=\ln(t)$ if $r= s$, we can see that $A_{f,p}=G_{r,s;\lambda}$. Using \lem{DB}, it follows that
\Eq{*}{
  \partial_\ell G_{r,s;\lambda}^\Delta(t)=\lambda_\ell  \qquad(\ell\in\{1,\dots,n\},\,t\in\R_+),
}
and 
\Eq{Glm}{
  \partial_\ell\partial_m G_{r,s;\lambda}^\Delta(t)
  =\lambda_m(\delta_{\ell,m}-\lambda_\ell)\frac{r+s-1}{t}\qquad(\ell,m\in\{1,\dots,n\},\,t\in\R_+).
}
Furthermore,
\Eq{Gyu}{
  \frac{p_0(y)(f(y)-f(u))}{p_0(u)f'(u)}
  =u\,\chi_{r,s}\Big(\frac{y}{u}\Big)\qquad(u,y\in\R_+).
}

\Thm{G1}{Let $n,k\in\N$ with $n,k\geq 2$, $\lambda\in\R_+^n$, $I_1,\dots,I_k$ be nonempty open subintervals of $\R_+$, $I:=I_1\times\dots\times I_k$, $(r_0,s_0),\dots,(r_k,s_k)\in\R^2$ and $\Phi\colon I\to\R_+$ be twice differentiable with nonvanishing first derivatives.
Then, for the inequality
\Eq{GHM}{
  G_{r_0,s_0;\lambda}(\Phi(x_1),\dots,\Phi(x_n))
  \leq \Phi\big(G_{r_1,s_1;\lambda}(x^1),\dots,G_{r_k,s_k;\lambda}(x^k)\big)
}
to be valid in the local sense it is necessary that
the $k\times k$ matrix $\Gamma(y)$ given by 
\Eq{*}{
    \Gamma(y):=\bigg(-\partial_i\partial_j\Phi(y)
    -\partial_j\Phi(y)\partial_i\Phi(y)\frac{r_0+s_0-1}{\Phi(y)}
    +\delta_{i,j}\partial_j\Phi(y)\frac{r_j+s_j-1}{y_j}\bigg)_{i,j=1}^k
    }
be positive semidefinite for all $y\in I$. Conversely, if this matrix is positive definite for all $y\in I$, then \eq{GHM} holds in the local sense on $I$.
}

\begin{proof}
	The necessity of the positive semidefiniteness of $\Gamma(y)$ for all $y\in I$ is an immediate consequence of \thm{17PZloc}, \thm{NecC}, and formula \eq{Glm}. The other direction is also obvious due to \thm{SP}.
\end{proof}

\Thm{Ggl}{Let $k\in\N$ with $k\geq2$, $I_1,\dots,I_k$ be nonempty open subintervals of $\R_+$, $I:=I_1\times\dots\times I_k$, $(r_0,s_0),\dots,(r_k,s_k)\in\R^2$  and let $\Phi\colon I\to\R_+$ be partially differentiable. Assume that, for all $(u,y)\in I^2$, the inequality 
\Eq{GSC0G}{
  \Phi(u)\chi_{r_0,s_0}\bigg(\frac{\Phi(y)}{\Phi(u)}\bigg)
  \leq \sum_{j=1}^k\partial_j\Phi(u)
  u_j\chi_{r_j,s_j}\Big(\frac{y_j}{u_j}\Big)
}
is valid. Then, for all $n\in\N$ and $\lambda\in\R_+^n$ with $\lambda_1+\dots+\lambda_n=1$, the inequality  \eq{GHM}
holds in the global sense on $I$.}

\begin{proof}
The statement directly follows from \thm{17PZgl} because inequality \eq{GSC0} turns out to be equivalent to \eq{GSC0G} by applying formula \eq{Gyu}.  
\end{proof}

For the investigation of the particular cases when $\Phi$ is the sum and the product function, we will need the following auxiliary result.

\Lem{PD}{Let $k\in\N$ with $k\geq 2$ and $c_0,c_1,\dots,c_k\in\R$. Then the matrix
\Eq{*}{
  C:=\big(\delta_{i,j}c_i+c_0\big)_{i,j=1}^k
}
is positive semidefinite if and only if either $c_0,c_1,\dots,c_k\geq0$ or there exists $i\in\{0,\dots,k\}$ such that $c_i<0$ and $c_j>0$ for all $j\in\{0,\dots,k\}\setminus\{i\}$ and 
\Eq{rec}{
  \frac{1}{c_0}+\frac{1}{c_1}+\cdots+\frac{1}{c_k}\leq0.
}
Furthermore, $C$ is positive definite if and only if either $c_0,c_1,\dots,c_k\geq0$ and $c_i=0$ can hold for at most one index $i\in\{0,\dots,k\}$ or there exists $i\in\{0,\dots,k\}$ such that $c_i<0$ and $c_j>0$ for all $j\in\{0,\dots,k\}\setminus\{i\}$ and \eq{rec} is valid with a strict inequality.}

\begin{proof} The quadratic form $Q\colon \R^k\to\R$ generated by $C$ is given by
\Eq{*}{
  Q(x):=c_0\Bigg(\sum_{\ell=1}^k x_\ell\Bigg)^2+\sum_{\ell=1}^k c_\ell x_\ell^2 \qquad (x\in\R^k).
}

Assume first that $Q$ is positive semidefinite. This  means that $Q(x)\geq0$ for all $x\in\R^k$. With the notation $x_0:=-\sum_{\ell=1}^kx_\ell$, this inequality can be rewritten as
\Eq{Q0}{
  \sum_{\ell=0}^k c_\ell x_\ell^2\geq0 \qquad\mbox{ for all } (x_0,x_1,\dots,x_k)\in\R^{k+1} \mbox{ with } x_0+x_1+\dots+x_k=0.
}
To prove the necessity of the condition, assume that $\min(c_0,c_1,\dots,c_k)<0$. Choose $i\in\{0,\dots,k\}$ such that $c_i=\min(c_0,c_1,\dots,c_k)<0$. For a fixed $j\in\{0,\dots,k\}\setminus\{i\}$, define the vector $x=(x_0,x_1,\dots,x_k)$ by $x_\ell:=\delta_{i,\ell}-\delta_{j,\ell}$ for $\ell\in\{0,\dots,k\}$, that is, $x_\ell:=0$ if $\ell\in\{0,\dots,k\}\setminus\{i,j\}$ and $x_i:=1$, $x_j:=-1$. Thus $x_0+x_1+\dots+x_k=0$ holds, which, by \eq{Q0}, yields that $\sum_{\ell=0}^k c_\ell x_\ell^2\geq0$, that is, $c_i+c_j\geq0$. This implies that $c_j>0$ for all $j\in\{0,\dots,k\}\setminus\{i\}$.

To show that the inequality \eq{rec} is also valid, we substitute the vector $x=(x_0,x_1,\dots,x_k)$ given by
\Eq{1/c}{
  x_j:=\frac{1}{c_j} \qquad (j\in\{0,\dots,k\}\setminus\{i\})\qquad\mbox{and}\qquad x_i:=-\sum_{j\in\{0,\dots,k\}\setminus\{i\}}\frac{1}{c_j}.
}
Then, obviously, $x_0+x_1+\dots+x_k=0$, which, again by \eq{Q0}, yields $\sum_{\ell=0}^k c_\ell x_\ell^2\geq0$, that is,
\Eq{ci1}{
  \sum_{j\in\{0,\dots,k\}\setminus\{i\}}c_j\bigg(\frac{1}{c_j}\bigg)^2
  +c_i\Bigg(-\sum_{j\in\{0,\dots,k\}\setminus\{i\}}\frac{1}{c_j}\Bigg)^2\geq0.
}
After simplifications, this implies that
\Eq{ci}{
  1+c_i\Bigg(\sum_{j\in\{0,\dots,k\}\setminus\{i\}}\frac{1}{c_j}\Bigg)\geq0,
}
which, using that $c_i<0$, shows that the inequality \eq{rec} is valid.

Assume that $Q$ is positive definite. Then $Q(x)>0$ for all $x\in\R^k\setminus\{0\}$. With the notation $x_0:=-\sum_{\ell=1}^kx_\ell$ this inequality can be rewritten as
\Eq{Q+}{
  \sum_{\ell=0}^k c_\ell x_\ell^2>0 \qquad\mbox{ for all } (x_0,x_1,\dots,x_k)\in\R^{k+1}\setminus\{0\} \mbox{ with } x_0+x_1+\dots+x_k=0.
}
Since the positive definiteness of $Q$ implies its positive semidefiniteness, there are two possible cases:
\begin{enumerate}[(a)]
	\item $c_0,c_1,\dots,c_k\geq0$;
	\item there exists $i\in\{0,\dots,k\}$ such that $c_i<0$ and $c_j>0$ for all $j\in\{0,\dots,k\}\setminus\{i\}$ and \eq{rec} holds.
\end{enumerate}

Assume first that case (a) holds. If $c_i=c_j=0$ were valid for some $i,j\in\{0,\dots,k\}$ with $i\neq j$, then substituting $x_\ell:=\delta_{i,\ell}-\delta_{j,\ell}$ for $\ell\in\{0,\dots,k\}$ into \eq{Q+}, we would get that $0<\sum_{\ell=0}^k c_\ell x_\ell^2=c_i+c_j=0$. This contradiction shows that $c_i=0$ can hold only for at most one index $i\in\{0,\dots,k\}$. 

Consider now the case (b). Substituting the vector $x=(x_0,x_1,\dots,x_k)$ given by \eq{1/c} into \eq{Q+}, it follows that $\sum_{\ell=0}^k c_\ell x_\ell^2>0$. This implies that the inequalities \eq{ci1} and then \eq{ci} hold with strict inequalities. Therefore, \eq{rec} is also satisfied with a strict inequality.

Now we show the sufficiency of the conditions. In case (a), it is clear that the inequality \eq{Q0} holds, whence $Q$ is positive semidefinite. If, in addition, $c_i=0$ holds for at most one index $i\in\{0,\dots,k\}$, then \eq{Q+} is also valid, that is, $Q$ is positive definite in this case.

In case (b), it follows from \eq{rec} that 
\Eq{ci2}{
    c_i\geq -\Bigg(\sum_{j\in\{0,\dots,k\}\setminus\{i\}}\frac{1}{c_j}\Bigg)^{-1}.
}
Let $(x_0,x_1,\dots,x_k)\in\R^{k+1}$ with $x_0+x_1+\dots+x_k=0$. Applying the Cauchy--Schwarz inequality to the vectors 
\Eq{*}{
   \bigg(\frac{1}{\sqrt{c_j}}\bigg)_{j\in\{0,\dots,k\}\setminus\{i\}}
   \qquad\mbox{and}\qquad
  \Big(-\sqrt{c_j}x_j\Big)_{j\in\{0,\dots,k\}\setminus\{i\}},
}
we obtain that
\Eq{*}{
  \Bigg(\sum_{j\in\{0,\dots,k\}\setminus\{i\}}\frac{1}{c_j}\Bigg)\Bigg(\sum_{j\in\{0,\dots,k\}\setminus\{i\}}c_jx_j^2\Bigg)
  \geq \Bigg(\sum_{j\in\{0,\dots,k\}\setminus\{i\}}-x_j\Bigg)^2=x_i^2.
}
Therefore, combining this inequality with \eq{ci2}, we can conclude that
\Eq{cii}{
  \sum_{j=0}^kc_jx_j^2\geq
  \sum_{j\in\{0,\dots,k\}\setminus\{i\}}c_jx_j^2
  -\Bigg(\sum_{j\in\{0,\dots,k\}\setminus\{i\}}\frac{1}{c_j}\Bigg)^{-1}x_i^2\geq0.
}
Thus, we have proved that \eq{Q0} is valid and hence $Q$ is positive semidefinite. If in this case \eq{rec} is valid with a strict inequality, then \eq{ci2} is also strict. Then, for a nonzero vector $(x_0,x_1,\dots,x_k)\in\R^{k+1}$ with $x_0+x_1+\dots+x_k=0$ the first inequality in \eq{cii} is strict provided that $x_i\neq0$ and the second inequality is strict if $x_i=0$. Thus, \eq{Q+} is valid, which shows that $Q$ is positive definite.
\end{proof}

\subsection{Minkowski-type inequalities}
Our next result characterizes the Minkowski-type inequality for Gini means in the local sense.

\Thm{G1+}{Let $n,k\in\N$ with $n,k\geq 2$, $\lambda\in\R_+^n$, $I_1,\dots,I_k$ be nonempty open subintervals of $\R_+$, $I:=I_1\times\dots\times I_k$, $(r_0,s_0),\dots,(r_k,s_k)\in\R^2$, $\gamma_i:=r_i+s_i-1$ for $i\in\{0,\dots,k\}$. For the inequality
\Eq{GM}{
  G_{r_0,s_0;\lambda}(x_1^1+\dots+x_1^k,\dots,x_n^1+\dots+x_n^k)
  \leq G_{r_1,s_1;\lambda}(x^1)+\dots+G_{r_k,s_k;\lambda}(x^k)
}
to hold in the local sense on $I$, it is necessary that exactly one of the following cases be valid:
\begin{enumerate}[(i)]
	\item \Eq{GM0}{
		\gamma_0\leq 0\leq\min(\gamma_1,\dots,\gamma_k);
	}
	\item $\gamma_0,\gamma_1,\dots,\gamma_k>0$ and
	\Eq{G1+}{
		\sum_{i\in J_+}\Big(\frac{1}{\gamma_i}-\frac{1}{\gamma_0}\Big)\sup I_i\leq \sum_{i\in J_-} \Big(\frac{1}{\gamma_0}-\frac{1}{\gamma_i}\Big)\inf I_i;
	}
	\item $\gamma_0<0$ and there exists $i\in\{1,\dots,k\}$ such that $\gamma_i<0$, for all $j\in\{1,\dots,k\}\setminus\{i\}$, $\gamma_j>0$, and inequality \eq{G1+} is also valid,
\end{enumerate} 
where, for the last two cases, we define
\Eq{*}{
		J_+:=\bigg\{i\in\{1,\dots,k\}\colon \frac1{\gamma_i}>\frac{1}{\gamma_0}\bigg\} \qquad \mbox{and}\qquad J_-:=\bigg\{i\in\{1,\dots,k\}\colon \frac1{\gamma_0}>\frac{1}{\gamma_i}\bigg\}.
}
Conversely, if either \eq{GM0} is valid and $\gamma_i=0$ can hold for at most one $i\in\{0,\dots,k\}$ or $\gamma_0\neq\gamma_\ell$ for some $\ell\in\{1,\dots,k\}$ and one of the conditions (ii) or (iii) hold, then \eq{GM} is valid in the local sense on $I$.}

\begin{proof} We apply \thm{G1} with the setting $I_0:=\R_+$ and $\Phi\colon I\to I_0$ defined by $\Phi(y):=y_1+\dots+y_k$. According to \thm{G1}, for the validity of \eq{GM} in the local sense, it is necessary that the values of the function $\Gamma\colon I\to\R^{k\times k}$ defined by
\Eq{*}{
	\Gamma(y):=\bigg(\delta_{i,j}\frac{\gamma_i}{y_i}-\frac{\gamma_0}{y_1+\dots+y_k}\bigg)_{i,j=1}^k
}
be positive semidefinite matrices. By \lem{PD}, this property is characterized by the following system of conditions: either 
\Eq{c0i}{
	c_0:=-\frac{\gamma_0}{y_1+\cdots+y_k}\geq0, \qquad c_i:=\frac{\gamma_i}{y_i}\geq0 \qquad (i\in\{1,\dots,k\});
}
or there exists $i\in\{0,\dots,k\}$ such that $c_i<0$, for all $j\in\{0,\dots,k\}\setminus\{i\}$, $c_j>0$, and \eq{rec} holds, i.e., 
\Eq{*}{
	\frac{1}{c_0}+\frac{1}{c_1}+\cdots+\frac{1}{c_k}\leq0.
}
Observe that $\sign(\gamma_0)=-\sign(c_0)$ and $\sign(\gamma_i)=\sign(c_i)$ for all $i\in\{1,\dots,k\}$.
Therefore, the first alternative can hold if and only if \eq{GM0} is satisfied. The second alternative can be valid if and only if either $\gamma_0,\gamma_1,\dots,\gamma_k>0$ or $\gamma_0<0$ and there exists $i\in\{1,\dots,k\}$ such that $\gamma_i<0$, for all $j\in\{1,\dots,k\}\setminus\{i\}$, $\gamma_j>0$, and 
\Eq{*}{
  0\leq \Big(\frac{1}{\gamma_0}-\frac{1}{\gamma_1}\Big)y_1
  +\dots+ \Big(\frac{1}{\gamma_0}-\frac{1}{\gamma_k}\Big)y_k
  \qquad (y\in I).
}
One can easily see that this inequality can be rewritten as \eq{G1+} and hence the necessity of the other two alternatives has been established.

To prove the reverse implication of the theorem, consider first the case when \eq{GM0} is valid and $\gamma_i=0$ for at most one $i\in\{0,\dots,k\}$. Then, for every $y\in I$, the numbers $c_0,c_1,\dots,c_k$ defined in \eq{c0i} are nonnegative and $c_i=0$ can hold for at most one $i\in\{0,\dots,k\}$. Thus, in view of the second assertion of \lem{PD}, it follows that $\Gamma(y)$ is positive definite.

Now consider the second case when, for some $\ell\in\{1,\dots,k\}$, we have that $\gamma_0\neq\gamma_\ell$ and either $\gamma_0,\gamma_1,\dots,\gamma_k>0$ or $\gamma_0<0$ and there exists $i\in\{1,\dots,k\}$ such that $\gamma_i<0$, for all $j\in\{1,\dots,k\}\setminus\{i\}$, $\gamma_j>0$ and \eq{G1+} is also valid. Let $y\in I$ be fixed. If $\ell\in J_+$, i.e., $\frac{1}{\gamma_\ell}>\frac{1}{\gamma_0}$, then 
\Eq{*}{
  \Big(\frac{1}{\gamma_\ell}-\frac{1}{\gamma_0}\Big)y_\ell
  <\Big(\frac{1}{\gamma_\ell}-\frac{1}{\gamma_0}\Big)\sup I_\ell,
}
while if $\ell\in J_-$, i.e., $\frac{1}{\gamma_0}>\frac{1}{\gamma_\ell}$, then 
\Eq{*}{
  \Big(\frac{1}{\gamma_0}-\frac{1}{\gamma_\ell}\Big)\inf I_\ell
  <\Big(\frac{1}{\gamma_0}-\frac{1}{\gamma_\ell}\Big)y_\ell.
}
Therefore, \eq{G1+} implies that
\Eq{*}{
  \sum_{i\in J_+}\Big(\frac{1}{\gamma_i}-\frac{1}{\gamma_0}\Big)y_i < \sum_{i\in J_-} \Big(\frac{1}{\gamma_0}-\frac{1}{\gamma_i}\Big)y_i,
}
which then yields that 
\Eq{*}{
  0<\Big(\frac{1}{\gamma_0}-\frac{1}{\gamma_1}\Big)y_1
  +\dots+ \Big(\frac{1}{\gamma_0}-\frac{1}{\gamma_k}\Big)y_k.
}
Hence, \eq{rec} is valid with a strict inequality sign. On the other hand, with the exception of one index, the numbers $c_0,\dots,c_k$ are positive. Thus, in view of the second assertion of \lem{PD}, it follows that $\Gamma(y)$ is positive definite in this case as well.
\end{proof}

\Cor{G1++}{Let $n,k\in\N$ with $n,k\geq 2$, $\lambda\in\R_+^n$, $I_1,\dots,I_k$ be nonempty open subintervals of $\R_+$ with $\inf I_1=\dots=\inf I_k=0$, $I:=I_1\times\dots\times I_k$, $(r_0,s_0),\dots,(r_k,s_k)\in\R^2$. In order that the inequality \eq{GM} be valid in the local sense on $I$, it is necessary that
\Eq{G1++}{
  \max(1,r_0+s_0)\leq\min(r_1+s_1,\dots,r_k+s_k).
}
Conversely, if this inequality is strict, then \eq{GM} holds in the local sense on $I$.}

\begin{proof}
	The proof is based on \thm{G1+}. Denote $\gamma_i:=r_i+s_i-1$ for $i\in\{0,\dots,k\}$. If condition \eq{GM0} holds, then \eq{G1++} is obvious because 
	\Eq{*}{
	  \max(1,r_0+s_0)=1+\max(0,\gamma_0)=1\leq 1+\min(\gamma_1,\dots,\gamma_k)=\min(r_1+s_1,\dots,r_k+s_k).
	}
	
	In the remaining two cases \eq{G1+} is valid. However, due to our assumptions on the intervals, the right hand side of \eq{G1+} is equal to 0. Therefore, the left hand side of this equality must be an empty sum, i.e., $J_+=\emptyset$, which means
	\Eq{g0i}{
		\frac{1}{\gamma_0}\geq\frac{1}{\gamma_i}
		\qquad(i\in\{1,\dots,k\})
	} should be valid. In the case $\gamma_0,\gamma_1,\dots,\gamma_k>0$, that is, when $1\leq\min(r_0+s_0,r_1+s_1,\dots,r_k+s_k)$, the inequalities in \eq{g0i} hold if and only if \eq{G1++} is satisfied. In the case when $\gamma_0<0$, then for at least one $j\in\{1,\dots,k\}$, we have that $\gamma_j>0$, and hence \eq{g0i} cannot hold.
	
    Assume now that \eq{G1++} is satisfied with a strict inequality. Then $\gamma_1,\dots,\gamma_k>0$. If $\gamma_0\leq 0$, then \eq{G1++} implies that \eq{GM0} is valid with a strict inequality and thus the first alternative of the sufficiency of \thm{G1+} holds. If $\gamma_0>0$, then the strict version of \eq{G1++} shows that $\gamma_0<\gamma_i$ for all $i\in\{1,\dots,k\}$, therefore, the left and the right hand sides of \eq{G1+} are equal to zero and the second alternative of the sufficiency of \thm{G1+} holds. The third alternative of the sufficiency of \thm{G1+} cannot happen if \eq{G1++} is valid.
\end{proof}

For the global validity of the Minkowski-type inequality, \thm{Ggl} establishes the following sufficient condition.

\Thm{GglM}{Let $k\in\N$ with $k\geq2$, $I_1,\dots,I_k$ be nonempty open subintervals of $\R_+$, $I:=I_1\times\dots\times I_k$, $(r_0,s_0),\dots,(r_k,s_k)\in\R^2$. Assume that, for all $(u,y)\in I^2$, the inequality 
\Eq{GSC0GM}{
  \chi_{r_0,s_0}\Big(\frac{y_1+\dots+y_k}{u_1+\dots+u_k}\Big)
  \leq \sum_{j=1}^k
  \frac{u_j}{u_1+\dots+u_k}\chi_{r_j,s_j}\Big(\frac{y_j}{u_j}\Big)
}
holds. Then, for all $n\in\N$ and $\lambda\in\R_+^n$ with $\lambda_1+\dots+\lambda_n=1$, the inequality \eq{GM}
holds in the global sense on $I$.}

\begin{proof}
With $\Phi(y_1,\dots,y_k):=y_1+\dots+y_k$, the condition \eq{GSC0G} turns out to be equivalent to \eq{GSC0GM} and hence the result follows from \thm{Ggl}.
\end{proof}

\Cor{GglM}{Let $k\in\N$ with $k\geq2$ and $(r_0,s_0),\dots,(r_k,s_k)\in\R^2$. Assume that, for all $z\in\R_+^k$ and $t_1,\dots,t_k\in[0,1]$ with $t_1+\dots+t_k=1$, the following inequality is valid 
\Eq{GSC0GM+}{
  \chi_{r_0,s_0}(t_1z_1+\dots+t_kz_k)
  \leq \sum_{j=1}^k t_j\chi_{r_j,s_j}(z_j).
}
Then, for all $n\in\N$ and $\lambda\in\R_+^n$ with $\lambda_1+\dots+\lambda_n=1$, the inequality \eq{GM}
holds in the global sense on $\R_+^k$.}

\begin{proof}
Let $(u,y)\in (\R_+^k)^2$ be arbitrary. Then, with the substitutions
\Eq{*}{
  z_j:=\frac{y_j}{u_j}\qquad \mbox{and}\qquad
  t_j:=\frac{u_j}{u_1+\dots+u_k} \qquad(j\in\{1,\dots,k\})
}
inequality \eq{GSC0GM+} implies \eq{GSC0GM}. Therefore, \eq{GSC0GM} holds for all $(u,y)\in (\R_+^k)^2$ and, according to \thm{GglM}, this condition yields that the inequality  \eq{GM} holds in the global sense on $\R_+^k$.
\end{proof}

In order to compare our results above to existing ones, we recall two theorems related to the global validity of the Minkowski-type inequalities. In the setting of two-variable Gini means the Minkowski inequality was characterized by Czinder and Páles in \cite[Theorem 5]{CziPal00} (see also \cite{LosPal97a} for a particular case).

\Thm{CP}{
	Let $k\in\N$ with $k\geq 2$, $(r_0,s_0),\dots,(r_k,s_k)\in\R^2$. Then the inequality 
	\Eq{M2vG}{
		G_{r_0,s_0}(x_1+\dots+x_k,y_1+\dots+y_k)
		\leq G_{r_1,s_1}(x_1,y_1)+\dots+G_{r_k,s_k}(x_k,y_k)
	}
	is valid for all $x_1,\dots,x_k,y_1,\dots,y_k\in\R_+$ if and only if 
	\begin{enumerate}[(i)]
      \item $0\leq\min(r_1,s_1,\dots.r_k,s_k)$,
      \item $\min(r_0,s_0)\leq \min(1,r_1,s_1,\dots.r_k,s_k)$,
      \item $\max(1,r_0+s_0)\leq\min(r_1+s_1,\dots,r_k+s_k)$.
	\end{enumerate}
}

\Rem{mink}{
Observe that the third condition in the above theorem is the necessary condition for the local validity of \eq{M2vG} on $\R_+$. The conditions (i) and (ii) are, however, not necessary for the local validity of \eq{M2vG} on $\R_+$.
}

Necessary and sufficient conditions for the global validity of the Minkowski-type inequality for Gini means with arbitrary number of variables was established by Páles in \cite[Theorem 3.1]{Pal82b}.

\Thm{Pal82}{
	Let $k\in\N$ with $k\geq 2$, $(r_0,s_0),\dots,(r_k,s_k)\in\R^2$. Then the inequality 
    \Eq{MnvG}{
    G_{r_0,s_0}(x_1^1+\dots+x_1^k,\dots,x_n^1+\dots+x_n^k)
    \leq G_{r_1,s_1}(x^1)+\dots+G_{r_k,s_k}(x^k)
    }
	is valid for all $n\in\N$ and $x\in\R_+^{n\times k}$ if and only if 
	\begin{enumerate}[(i)]
      \item $0\leq\min(r_1,s_1,\dots,r_k,s_k)$,
      \item $\min(r_0,s_0)\leq \min(1,r_1,s_1,\dots,r_k,s_k)$,
      \item $\max(1,r_0,s_0)\leq\min(\max(r_1,s_1),\dots,\max(r_k,s_k))$.
	\end{enumerate}
}

\Rem{mink2}{
Observe that conditions (i) and (ii) of the above two theorems are identical, therefore they may be necessary for the global validity of \eq{MnvG} for any fixed $n\in\N$. The form of the third condition related to any fixed $n\in\N$ is not known. We also note that conditions (i)-(iii) of \thm{Pal82} are also necessary and sufficient for the validity of the inequality \eq{GSC0GM+} on the domain indicated in \cor{GglM}.
}

\subsection{Hölder-type inequalities} Our next results characterize Hölder-type inequalities for Gini means in the local and in the global sense.

\Thm{GH}{Let $n,k\in\N$ with $n,k\geq 2$, $\lambda\in\R_+^n$, $I_1,\dots,I_k$ be nonempty open subintervals of $\R_+$, $I:=I_1\times\dots\times I_k$, $(r_0,s_0),\dots,(r_k,s_k)\in\R^2$, $\gamma_i:=r_i+s_i$ for $i\in\{0,\dots,k\}$. Then, in order that the inequality
\Eq{GH}{
  G_{-r_0,-s_0;\lambda}(x_1^1\cdots x_1^k,\dots,x_n^1 \cdots x_n^k)
  \leq G_{r_1,s_1;\lambda}(x^1)\cdots G_{r_k,s_k;\lambda}(x^k)
}
be valid in the local sense on $\R_+^k$ it is necessary that 
either $\gamma_0,\gamma_1,\dots,\gamma_k\geq0$ 
or there exists $i\in\{0,\dots,k\}$ such that $\gamma_i<0$, for all $j\in\{0,\dots,k\}\setminus\{i\}$, $\gamma_j>0$ and 
\Eq{recg}{
  \frac{1}{\gamma_0}+\frac{1}{\gamma_1}+\dots+\frac{1}{\gamma_k}\leq0
}
holds. Conversely, if either $\gamma_0,\gamma_1,\dots,\gamma_k\geq0$ and $\gamma_i=0$ for at most one index $i\in\{0,\dots,k\}$ or there exists $i\in\{0,\dots,k\}$ such that $\gamma_i<0$, for all $j\in\{1,\dots,k\}\setminus\{i\}$, $\gamma_j>0$ and \eq{recg} is valid with a strict inequality, then \eq{GH} holds in the local sense on $\R_+^k$.}

\begin{proof} We apply \thm{G1} with the function $\Phi\colon I\to\R_+$ defined by $\Phi(y):=y_1\cdots y_k$. 
Then for the validity of \eq{GH} in the local sense it is necessary (and sufficient) that the values of the function $\Gamma\colon \R_+^k\to\R^{k\times k}$ defined by
\Eq{*}{
    \Gamma(y)&:=\bigg((\delta_{i,j}-1)
    \frac{1}{y_iy_j}\prod_{\ell=1}^ky_\ell
    -(-\gamma_0-1)\frac{1}{y_iy_j}\prod_{\ell=1}^ky_\ell 
    +\delta_{i,j}(\gamma_j-1)\frac{1}{y_iy_j}\prod_{\ell=1}^ky_\ell\bigg)_{i,j=1}^k\\
    &\phantom{:}=\bigg(\frac{1}{y_iy_j}\prod_{\ell=1}^ky_\ell\big(\delta_{i,j}\gamma_j+\gamma_0\big)\bigg)_{i,j=1}^k
    }
be positive semidefinite (positive definite) matrices for all $y\in I$. However, this property holds if and only if the scalar matrix
\Eq{*}{
    \Gamma^*&:=\big(\delta_{i,j}\gamma_j+\gamma_0\big)_{i,j=1}^k
}
is positive semidefinite (positive definite). The statement now follows from \lem{PD} with $c_i:=\gamma_i$ ($i\in\{0,\dots,k\}$).
\end{proof}

\Thm{GglH}{
	Let $k\in\N$, $k\geq2, I_1,\dots,I_k\subseteq\R_+$ be nonempty open intervals, $(r_0,s_0),\dots,(r_k,s_k)\in\R^2$. Assume that, for all $z_1\in (I_1/I_1),\dots,z_k\in (I_k/I_k)$, the inequality
\Eq{GSC0GH}{
	\chi_{-r_0,-s_0}(z_1\cdots z_k)\leq\sum_{j=1}^{k}\chi_{r_j,s_j}(z_j)
}
	holds. Then, for all $n\in\N$ and $\lambda\in\R_+^n$ with $\lambda_1+\dots+\lambda_n=1$, the inequality \eq{GH} holds in the global sense on $I$. 
}

\begin{proof}
	With the function $\Phi(y_1,\dots,y_k):=y_1\cdots y_k$, condition \eq{GSC0G} turns out to be equivalent to 
	\Eq{*}{
	\chi_{-r_0,-s_0}\Big(\frac{y_1\cdots y_k}{u_1\cdots u_k}\Big)\leq\sum_{j=1}^{k}\chi_{r_j,s_j}\Big(\frac{y_j}{u_j}\Big).
}
Introducing the new variables $z_i:=y_i/u_i$ for $i\in\{1,\dots,k\}$, we can conclude that \eq{GSC0GH} is valid for all $z_1\in (I_1/I_1), \dots, z_k\in (I_k/I_k)$ if and only if the above inequality holds for all $(y,u)\in (I_1\times\dots\times I_k)^2$. Hence the result follows from \thm{Ggl}.
\end{proof}

The global validity of \eq{GH} with a non-fixed number of variables was characterized by Páles in \cite{Pal83a,Pal83c}. 

\Thm{GHg}{
	Let $k\in\N$ with $k\geq 2$, $(r_0,s_0),\dots,(r_k,s_k)\in\R^2$. Then the inequality 
	\Eq{*}{
		G_{-r_0,-s_0}(x_1^1\cdots x_1^k,\dots,x_n^1 \cdots x_n^k)
		\leq G_{r_1,s_1}(x^1)\cdots G_{r_k,s_k}(x^k)
	}
	is valid for all $n\in\N$ and $x\in\R_+^{n\times k}$ if and only if
	\begin{enumerate}[(i)]
	 \item for all $i\in\{0,\dots,k\}$, $\max(s_i,r_i)\geq0$ and
	 \item for all $i\in\{0,\dots,k\}$ with $\min(s_i,r_i)<0$, we have $\max(s_j,r_j)>0$ for all $j\in\{0,\dots,k\}\setminus\{i\}$ and
	\Eq{*}{
     \frac{1}{\min(s_i,r_i)}+\sum_{\substack{j=0\\j\neq i}}^{k}\frac{1}{\max(s_j,r_j)}\leq 0.
	}
	\end{enumerate}
}

\Rem{h}{
We note that the conditions (i) and (ii) are necessary and sufficient for the validity of the inequality \eq{GSC0GH} for all $z_1,\dots,z_k\in\R_+$ (cf.\ \cite{Pal83c}).
}

\section*{Acknowledgement}

The authors wish to thank the anonymous referee for her/his detailed comments and suggestions. 

\def\MR#1{}

\end{document}